\documentclass[a4paper,12pt]{article}
\usepackage[dvips]{epsfig}
\usepackage{amsmath,amssymb,euscript,graphicx,amsfonts}
\usepackage{enumerate,color}
\usepackage{grap hicx}
\usepackage{hyperref, minitoc}
\newtheorem{theorem}{Theorem}[section]
\newtheorem{proposition}[theorem]{Proposition}
\newtheorem{lemma}[theorem]{Lemma}
\newtheorem{corollary}[theorem]{Corollary}

\newcommand{\OO}{\mathcal{O}}
\newcommand{\ZZ}{\mathbb{Z}}
\renewcommand{\S}{\mathcal{S}}
\renewcommand{\P}{\mathcal{P}}
\def\demo{{\bf Proof}\hskip10pt} \def\qqed{\hfill $\Box$}
\def\di{\bigm|} \def\lg{\langle} \def\rg{\rangle}
\def\PSL{\hbox{\rm PSL}} \def\a{\alpha} \def\f{\noindent} \def\FF{\mathbb{F}}
\def\Aut{\hbox{\rm Aut\,}} \def\PG{\hbox{\rm PG}} \def\th{\theta} \def\O{\Omega}

\begin{document}

\begin{center}
{\bf\large On Hamilton paths in vertex-transitive graphs of order $10p$
}
\footnote{This work was supported in part  by the National Natural Science Foundation of China
 (12071312,12471332).}
\end{center}


\begin{center}
Shaofei Du, Wenjuan Luo  and  Hao Yu$^{*}$\\
\medskip
{\it {\small
School of Mathematical Sciences \\
Capital Normal University\\
Bejing 100048, People's Republic of China\\
}}
\end{center}

\renewcommand{\thefootnote}{\empty}

\footnotetext{$^{\ast}$  Corresponding author; Email: 3485676673@qq.com}

\begin{abstract}
  It was shown by Kutnar, Maru\v si\v c and Zhang in 2012 that every connected vertex-transitive graph of order $10p$, where $p$ is a prime and $p\ne 7$, contains a Hamilton path, except for graphs $X$ arising from the action of $\PSL(2, s^m)$ on cosets of $\ZZ_s^m\rtimes \ZZ_{\frac{s^m-1}{10}}$, where $s$ is a prime. In this paper, Hamilton cycles of these exceptions $X$ will be  found.
\end{abstract}

\section{Introduction}\label{sec:intro}
\indent

A finite simple path (resp. cycle) going through all vertices of the graph is called {\it a Hamilton path (resp. cycle)}.
In 1969, Lov\'asz \cite{L70} asked if there exists a finite connected vertex-transitive graph without a Hamilton path;
and in 1981, Alspach \cite{A81} asked if there exists an infinite number of connected vertex-transitive graphs that do not have a hamiltonian cycle.

Till now, only four connected vertex-transitive graphs of order at least 3 that do not have a Hamilton cycle are known to exist.
These four graphs are the Petersen graph, the Coxeter graph and the two graphs obtained
from them by replacing each vertex by a triangle.
The fact that none of these four graphs is a Cayley graph  has led to a folklore conjecture that
every Cayley graph contains a Hamilton cycle (see  \cite{A89,D83,GWM14,GKMM12,GM07,KM09,
DM83,WM15,WM18} and the survey paper \cite{CG96} for the current status).

For vertex-transitive graphs, the existence of Hamilton paths, even Hamilton cycles, has been confirmed
for graphs of particular orders, such as, $kp$ where $k\le 6$, $p^j$ where $j\le 5$ and $2p^2$ (see \cite{C98,KM08,KMZ12,KS09,DM85,DM87,MP82,MP83,Z15} and the survey paper \cite{KM09}).
Recently, Kutnar, Maru\v si\v c and the first author proved that every connected vertex-transitive graph of order $pq$, where $p$ and $q$ are primes, has a Hamilton cycle,
except for the Peterson graph (see\cite{DKD1,DKD2}).
As for the vertex-transitive graphs of order $2pq$, where both $p$ and $q$ are primes, Tian, Yu and the
first author \cite{DTY} showed that every primitive graph of such order contains a Hamilton cycle, except for the Coxeter graph.
But   it is still open for the imprimitive case.
In \cite{KMZ12}, Kutnar, Maru\v si\v c and  Zhang showed that every connected vertex-transitive graph of
order $2\cdot 5\cdot p$  where $p$ is a prime
and $p\ne 7$, contains a Hamilton path, except for  graphs $X$ arising from the action of $\PSL(2,s^m)$
on cosets of $\ZZ_s^m\rtimes \ZZ_{\frac{s^m-1}{10}}$  where $s^m+1=2p$.
In this paper, these exceptions $X$ will be studied and Hamilton cycles of that will be found.

\begin{theorem}\label{the:main}
Let $X$ be a connected graph whose automorphism group contains  a vertex-transitive subgroup  $\PSL(2,s^m)$  where $s$ is a prime, having a point stabilizer $\ZZ_s^m\rtimes\ZZ_{\frac{s^m-1}{10}}$,
where $s^m+1=2p$ for a  prime $p$. Then $X$ contains a  Hamilton cycle.
\end{theorem}

Combining Theorem~\ref{the:main} and the main result in \cite{KMZ12} just mentioned above, we get

\begin{corollary}
Every connected vertex-transitive graph of order $10p$ contains a Hamilton path, where $p$ is a prime and $p\ne 7$.
\end{corollary}

After this introductory section, some notations, basic definitions and useful facts will
be given in Section 2 and Theorem \ref{the:main} will be proved in Section 3.

\section{Preliminaries}
\label{sec:pre}
\noindent

Throughout this paper graphs are finite and undirected.
By $p$ we always denote a prime.
Given a graph $X$, by $V(X),~E(X)$ and $\Aut(X)$ we denote the vertex set, the edge set and the automorphism group of $X$ ,respectively.
Let $U$ and $W$ be two disjoint subsets of $V(X)$.
By $X\lg U \rg$ and $X[U,W]$ we denote the subgraph of $X$ induced by $U$ and the bipartite subgraph with two biparts $U$ and $W$, respectively.
In the case when $X\lg U \rg$ and $X[U,W]$ are regular, $d(U)$ and $d(U,W)$ denote  the valency of $X\lg U\rg$ and $X[U,W]$, respectively.

In what follows, we recall some definitions and known facts.

\vskip 3mm
\f {\bf (1)  Generalized orbital graphs}

\vskip 3mm
A permutation group $G$ on a finite set $V$ induces a natural action of $G$ on $V\times V$, whose orbits are called {\em orbitals} of $G$.
Moreover, the orbital $\OO_0=\{(x,x) \colon x\in V\}$ is said to be {\it trivial.}
Every orbital $\OO_i$ has a paired orbital $\OO_i^*=\{ (y,x)\di (x,y)\in \OO_i\}$ and $\OO_i$ is said to be {\em self-paired} if it coincides with $\OO_i^*$.
Given an arbitrary union $\OO$ of some orbitals such that $\OO\cap \OO_0= \emptyset$, we get a simple digraph $X(G, \OO)$ which is called orbital digraph, with the vertex set $V$ and the edge set $\cal{O}$.
This graph may be viewed as an undirected graph provided that $\OO_i^*$ is contained in $\OO$ for every orbital $\OO_i\subseteq \OO$.
It is well-known that every vertex-transitive graph can be described as an orbital graph. The orbital graph $X(G,  \OO_i\cup \OO_i^*)$ is called {\it basic}, where  $\OO_i$ is a non-trivial single orbital (self-paired or non-self-paired).

Let $H=G_\a$ for some $\a\in V$.
Then there is an 1-1 correspondence between the orbits of the action of $H$ on $\O$, called {\em suborbits} of $G$, and orbitals.
A suborbit corresponding to a self-paired orbital is said to be
{\em self-paired}.

\vskip 3mm
\f {\bf (2)  Semiregular automorphisms and quotient (multi)graphs}

\vskip 3mm
Let $m\geq 1$ and $n\geq 2$ be integers. An automorphism $\rho$ of a graph $X$ is called $(m,n)$-{\em semiregular} (in short, {\em semiregular}) if as a permutation on $V(X)$ it has a cycle decomposition consisting of $m$ cycles of length $n$.
The question whether all vertex-transitive graphs admit a semiregular automorphism  is one of famous open problems in algebraic graph theory (see, for example, \cite{bcc15,seven,DMMN07,G1,M81}).
Let $\P$ be the set of orbits of $\lg \rho \rg $.
Let $X_\P$ be the {\em quotient graph corresponding to $\P$}, the graph  whose vertex set is $\P$, with $A, B \in \P$ adjacent if there exist adjacent vertices $a \in A$ and $b \in B$ in $X$.
Let $X_\rho$ be the {\em quotient multigraph corresponding to $\rho$},  the multigraph  whose vertex set is $\P$ and in which $A,B \in \P$ are joined by $d(A,B)$ edges.
Note that the quotient graph $X_\P$ is precisely the underlying graph of $X_\rho$.

\vskip 3mm
\f {\bf (3)  Lifting cycle technique}
\vskip 3mm

If $G$ is primitive on $\O$, then  every basic orbital graph is  connected.
Clearly,  to prove that every vertex-transitive graph arising from a primitive group $G$ on $\O$
contains a Hamilton cycle, it suffices to show that  every basic orbital graph contains a Hamilton cycle.

One of tools on Hamilton cycle problems is the so-called lifting cycle technique (see \cite{A89,KM09,DM83}).
When the quotient is applied relative to a semiregular automorphism of prime order and the corresponding quotient multigraph possesses two adjacent orbits linked by a double edge encompassed within a Hamilton cycle, lifts of Hamilton cycles from quotient graphs are invariably achievable.
This double edge enables us to conveniently ``change direction" to procure a walk in the quotient that elevates into a full cycle above.

Let $X$ be a graph that admits an $(m,n)$-semiregular automorphism $\rho$.
Let $\P = \{S_1, S_2, \ldots , S_m\}$
be the set of orbits of $\rho$, and let $\pi : X \to X_\P$
be the corresponding projection of $X$ to its quotient $X_\P$.
For a (possibly closed) path $W = S_{i_1}S_{i_2}\ldots S_{i_k}$ in $X_\P$ we let the {\em lift} of $W$ be the set of all paths in $X$ that project to $W$.
The proof of following lemma is straightforward and is just a reformulation of \cite[Lemma~5]{MP82}.

\begin{lemma} \label{lem:cyclelift}
Let $X$ be a graph admitting
an $(m,p)$-semiregular automorphism $\rho$, where $p$ is a prime.
Let $C$ be a cycle of length $k$ in the quotient graph $X_\P$,
where $\P$ is the set of orbits of $\rho$.
Then, the lift of $C$ either contains a cycle of length
$kp$ or it consists of $p$ disjoint $k$-cycles.
In the latter case we have $d(S,S') = 1$ for every edge $SS'$ of $C$.
\end{lemma}

\vskip 3mm
\f {\bf (4) A result  on number theory}
\vskip 3mm

By $\FF$ we denote a finite field of order $q=p^k$, for a prime $p$.
A {\it diagonal equation} over $\FF$  is an equation of the type
$$a_1 x_1^{k_1}+\cdots+a_n x_n^{k_n}=b$$
with positive integers $k_1, \cdots, k_n$, coefficients $a_1, \cdots, a_n \in \FF^*$, and $b \in \FF$.
For $n=2$, the following  result will be used later.

\begin{proposition}\label{num} {\rm \cite[Theorem 6.37]{LN1997}}
Let $N$ be the number  of solutions $(x_1, x_2)$ of the diagonal equation $a_1 x_1^{k_1}+a_2 x_2^{k_2}=b$, where  $a_1, a_2, b\in \FF_q^*$.
Then
$$|N-q| \leq[(d_1-1)(d_2-1)-(1-q^{-\frac12}) M(d_1, d_2)] q^{\frac12},$$
where $d_i={\rm gcd}(k_i, q-1)$ for $i=1,2$ and
$M(d_1, d_2)$ is the number of pairs $(j_1, j_2)\in \ZZ^2$ such that $1\le j_i\le d_i-1$ and
$\frac{j_1}{d_1}+\frac{j_2}{d_2}\in \ZZ$.
\end{proposition}

\section{Proof of Theorem~\ref{the:main}}\label{sec:PGL}
\indent

To prove  Theorem~\ref{the:main}, let $k=s^m$ where $s$ is a prime  such that $k+1=2p$ for some prime
$p$ and $10\di(k-1).$ Set $\FF_k^*=\lg \th \rg $.
Then $k\equiv1({\rm mod}\,4)$ so that $\th^{\frac{k-1}2}=-1$. Set $G=\PSL(2,k)$ and

$$\begin{array}{ll}
&\ell=\overline{\left[\begin{array}{cc} 0&-1\\1&0\end{array}\right]},\,
t=\overline{\left[\begin{array}{cc} \th&0\\0&\th^{-1}\end{array}\right]},\,
u=\overline{\left[\begin{array}{cc} 1&1\\0&1\end{array}\right]},\,
 s(a,b)=\overline{\left[ \begin{array}{cc} a&b\\b\th&a\end{array}\right]},\\ &\\
 &S=\lg s(a,b)\di a, b\in\FF_k,  a^2-b^2\th=1 \rg.\end{array}$$
\f Then  $${\rm o} (\ell)=2, \, {\rm o} (t)=\frac{k-1}2,\,  {\rm o} (u)=s,\,  t^{\ell}=t^{-1}\, \, {\rm and}\, \, S\cong \ZZ_{\frac{k+1}2}.$$

Let $K=\lg u, t\rg \cong \ZZ_s^m\rtimes\ZZ_{\frac{k-1}2}$. Then $K=G_{\infty}$, the point-stablizer in
$G$, relative to $\infty$ in the projective line $\PG(1, k)=\{\infty,\,0,\,1,\,2,\,\cdots,\,k-1\}$.
Clearly, $\PG(1, k)$ can be identified with $G\backslash K$, the set of right cosets of $G$ relative to $K$.
Let $H=\lg u, t^5\rg  \le K$ and $\Omega=G\backslash H$, the set of right cosets of $G$ relative to $H$.
Consider the action of $G$ on $\Omega$, where $|\Omega|=|G:H|=5(k+1)=10p$.
From now on, we let $\a=H\in \Omega$.
The following lemma determines the suborbits of $G$ relative to $\a$. Remind that for any subset $M$ of $G$, $\a^M$ denotes the set $\{\a^m\di m\in M\}$.

\begin{lemma}\label{orbit} Acting on $\Omega$,
\begin{enumerate}
  \item[\rm(1)] $G$ has five single point suborbits: $\{\a^{t^i}\}$ and five suborbits $\a^{t^i\ell H}$ of length $k$, where  $i\in \{0,1,2,3,4\}$.  Every suborbit $\a^{t^i\ell H}$ is self-paired;
  \item[\rm(2)] $S$ has ten orbits: $\a^{t^iS}$ and $\a^{t^i\ell S}$, where  $i\in \{0,1,2,3,4\}$. All of them are of length $\frac{k+1}2$.
\end{enumerate}
\end{lemma}
\demo (1)  Since the group $K$ has two orbits on $\Omega$ (as $G$ is 2-transitive on $\PG(1,k)$) and $K=\cup_{i=0}^{4}Ht^i$, we get that the group $H$ has five single point suborbits $\{\a^{t^i}\}$ and five suborbits $\a^{t^i\ell H}$ of length $k$, where $i\in \{0,1,2,3,4\}$.
Since $t^i \ell$ is an involution, every suborbit $\a^{t^i\ell H}$ is self-paired.

\vskip 3mm

(2) Note that the group $S$ acts semiregularly on $\PG(1, k)$.
Thus, $ K^g\cap S=1$ for any $g\in G$.
Since $H^g\cap S\leq K^g\cap S$ for any $g\in G$, the group $S$ acts semiregularly on $\Omega$.
Since the group $S$  has two orbits of equal length on $\PG(1, k)$, $\infty$ is not contained in the $S$-orbit containing 0 and $K=\cup_{i=0}^{4}Ht^i$,
we know  that the group $S$ has ten orbits on $\Omega$: $\a^{t^iS}$ and $\a^{t^i\ell S}$,  where
$i\in \{0,1,2,3,4\}$.
\qqed
\vskip 3mm
Note that there is an  1-1 correspondence between  suborbits $\a^{t^i\ell H}$ and orbitals $\OO_i=\{(\a,\a^{t^i\ell h})^g\di g\in G,h\in H\}$, where $i\in\{0,1,2,3,4\}$.
Before proving Theorem~\ref{the:main},  we first find a  Hamilton cycle for basic  orbital graphs $Y(i)=X(G,\OO_i)$ with $i\in\{0,1,2,3,4\}$.

\begin{lemma}\label{0}
For any $i\in\{0,1,2,3,4\}$, the basic  orbital graph $Y(i)=X(G,\OO_i)$ contains a Hamilton cycle.
\end{lemma}
\demo By Lemma~\ref{orbit}, the cyclic group $S$ acts semiregularly on $\Omega$ and has ten orbits on $\Omega$.
Let $Y(i)_\S$  be the quotient graph of $Y(i)$ induced by $S$, where $|S|$ is a prime $p=\frac{k+1}2$.
Then $X_\S$ has ten vertices. Once we may derive the following two facts:
 \vskip 3mm
 {\it  (i)}\, $Y(i)_\S$ is a complete graph;  and
 \vskip 3mm
 {\it  (ii)}\,  $d(A,  B)\ge 2$  for any two $S$-orbits $A$ and $B$,
\vskip 3mm
\f then by Proposition~\ref{lem:cyclelift}, every Hamilton cycle in $Y(i)_\S$ can be lifted to a Hamilton cycle of $Y(i)$. Moreover, to prove these two facts, it suffices to show
\begin{eqnarray}\label{Eq0} d(\a^{t^{n}S}, \a^{t^{j}S}), \,  d(\a^{t^{n}S}, \a^{t^{j}\ell S}), \,      d(\a^{t^{n}\ell S},\a^{t^{j}\ell S})\ge2,\end{eqnarray}
where $j, n\in\{0,1,2,3,4\}$.

Remind that the neighborhood of $\a$ is $\a^{t^i\ell H}=\{ \a^{t^i\ell h}\di h\in H\}$,
where $H=\lg u, t^5\rg$.
Then the neighborhood of $\a^{t^{n}}$ is
$$\begin{array}{lcl}
Y(i)_1(\a^{t^{n}})&=&\{\a^{t^i\ell ht^{n}}\di h\in H\}=\{\a^{t^{i-n}\ell h_1}\di h_1\in H\}\\
&=&\{\a^g\di g=\overline{\left[\begin{array}{cc} 0&-\th^{i-n}\\\th^{-(i-n)}&x\end{array}\right]}, x\in
\FF_k\}.\\
\end{array}$$
In what follows, we shall deal with three cases in Eq(\ref{Eq0}), separately.

\vskip 3mm
{\it Claim 1: $d(\a^{t^{n}S}, \a^{t^{j}S})\ge 2.$}
\vskip 3mm

Clearly,  $d(\a^{t^{n}S}, \a^{t^{j}S})=|Y(i)_1(\a^{t^{n}})\cap \a^{t^{j}S}|$, which is the number of
solutions $s(a,b)$  of the equation
\begin{eqnarray}\label{Eq1}
H\overline{\left[\begin{array}{cc} 0&-\th^{i-n}\\\th^{-(i-n)}&x\end{array}\right]}=Ht^js(a,b),
\end{eqnarray}
that is
$$\overline{\left[\begin{array}{cc} \th^{5r} &x_1\\ 0 &\th^{-5r} \end{array}\right]}
 \overline{\left[\begin{array}{cc} 0&-\th^{i-n}\\\th^{-(i-n)}&x\end{array}\right]}=t^js(a,b),$$
for some $r$ and $x_1$,  that is
$$\left[\begin{array}{cc} \th^{-(i-n)}x_1 &xx_1-\th^{5r+i-n}\\
 \th^{-(5r+i-n)}&\th^{-5r}x\end{array}\right]=\pm \left[ \begin{array}{cc}
 a\th^{j}&b\th^{j}\\b\th^{-j+1}&a\th^{-j}\end{array}\right].$$
Since $s(a,b)$ is a solution of  Eq(\ref{Eq1}) if and only if  $s(-a,-b)$ is a solution of  Eq(\ref{Eq1}), we just consider the  "+" case. Then
Eq(\ref{Eq1}) holds if and only if
$$\left\{\begin{array}{lcll} \th^{-(i-n)}x_1 &=&a\th^{j}, \, \hskip 2cm &(i)\\
xx_1-\th^{5r+i-n}&=&b\th^{j},\, &(ii) \\
 \th^{-(5r+i-n)}&=&b\th^{-j+1},\, &(iii) \\
 \th^{-5r}x&=&a\th^{-j},\, &(iv)\\
 a^2-b^2\th&=&1,  &(v)\end{array}\right.  $$
that is
$$\left\{\begin{array}{lcll}  x_1&=&a\th^{j+i-n},\,\hskip 2cm & (i')  \\
xx_1-\th^{5r+i-n}&=&b\th^{j},\, & (ii)\\
 b&=&\th^{-5r+j-i+n-1},\, &(iii')\\
x&=&a\th^{5r-j}, &(iv') \\
 a^2-b^2\th&=&1. &(v)\end{array}\right. $$
Inserting $(i')$, $(iii')$ and $(iv')$ to $(ii)$, we obtain
$$a^2\th^{5r+i-n}-\th^{5r+i-n}=\th^{-5r+2j-i+n-1}. \hskip 1cm (ii')$$
Note that $(ii')$ and $(iii')$ imply $(v)$.
Setting  $y=\th^{-r}$ and $c=-\th^{2(j-i+n)-1},$ the equation $(ii')$ becomes
\begin{eqnarray}\label{Eq2} a^2+cy^{10}=1.  \end{eqnarray}

Conversely, given any solution $(a,y)$ of Eq(\ref{Eq2}), we know that $(-a,y)$ is a solution too.
Then $x_1$, $b$ and $x$ are uniquely determined by $(i')$, $(iii')$ and $(iv')$, respectively; clearly,
$(ii')$ and so $(ii)$ holds; and finally $(v)$ holds.

In summary,  $d(\a^{t^{n}S}, \a^{t^{j}S})\ge 2$  if and only if  Eq(\ref{Eq2}) has solutions with $y\neq0$.
The remaining is to show the diagonal equation  Eq(\ref{Eq2}) has solutions by
using Proposition~\ref{num}.

Following the meaning of $M(d_1, d_2)$ in Proposition~\ref{num}, we have $d_1=2$, $d_2=10$ and $M(2, 10)$ is the number of pairs $(1, j_2)\in \ZZ^2$ such that $1\le j_2\le 9$ and $\frac{1}{2}+\frac{j_2}{10}\in \ZZ$.
Therefore, $M(2, 10)=1$. Let $N$ be the number of solutions of Eq(\ref{Eq2}) with $y\neq0$.
Since $(\pm1,0)$ are the only solutions with $y=0$, we see from Proposition~\ref{num} that
$$N\ge (k-8k^{\frac12}-1)-2=k-8k^{\frac12}-3.$$
Hence, Eq(\ref{Eq2}) has solutions  with $y\neq0$ if $k-8k^{\frac12}-3\ge 1$, which holds provided $k\ge 72$.
In the case $k\leq 71$, since $k+1=2p$ and $10\di k-1$, we get $k=61$.
Checking by Magma,  Eq(\ref{Eq2}) has solutions for any $c$, over $\FF_{61}$.

\vskip 3mm
{\it Claim 2: $d(\a^{t^{n}S}, \a^{t^{j}\ell S})\ge 2.$}
\vskip 3mm

Clearly, $d(\a^{t^{n}S}, \a^{t^{j}\ell S})=|Y(i)_1(\a^{t^{n}})\cap \a^{t^{j}\ell S}|$, which is
the number of solutions $s(a,b)$  of the equation
\begin{eqnarray}\label{Eq3}
H\overline{\left[\begin{array}{cc} 0&-\th^{i-n}\\\th^{-(i-n)}&x\end{array}\right]}=Ht^j\ell s(a,b),
\end{eqnarray}
that is (the left side of   Eq(\ref{Eq3}) is same with that in  Eq(\ref{Eq1}))
$$\left[\begin{array}{cc} \th^{-(i-n)}x_1 &xx_1-\th^{5r+i-n}\\
\th^{-(5r+i-n)}&\th^{-5r}x\end{array}\right]=\pm \left[ \begin{array}{cc}
-b\th^{j+1}&-a\th^{j}\\a\th^{-j}&b\th^{-j}\end{array}\right],$$
for some $r$ and $x_1$. Again we just need to consider the "+" case. So
Eq(\ref{Eq3}) holds if and only if
$$\left\{\begin{array}{lcll}
\th^{-(i-n)}x_1 &=&-b\th^{j+1}, \, \hskip 2cm &(i)\\
xx_1-\th^{5r+i-n}&=&-a\th^{j},\, &(ii) \\
 \th^{-(5r+i-n)}&=&a\th^{-j},\, &(iii) \\
 \th^{-5r}x&=&b\th^{-j},\, &(iv)\\
 a^2-b^2\th&=&1,  &(v)
 \end{array}\right.$$
that is
$$\left\{\begin{array}{lcll}
x_1&=&-b\th^{j+i-n+1},\,\hskip 2cm & (i')  \\
xx_1-\th^{5r+i-n}&=&-a\th^{j},\, & (ii)\\
 a&=&\th^{-5r+j-i+n},\, &(iii')\\
x&=&b\th^{5r-j}, &(iv') \\
 a^2-b^2\th&=&1. &(v)\end{array}\right. $$
Inserting $(i')$, $(iii')$ and $(iv')$ to $(ii)$, we obtain
$$b^2\th^{5r+i-n+1}+\th^{5r+i-n}=\th^{-5r+2j-i+n}. \hskip 1cm (ii')$$
Note that $(ii')$ and $(iii')$ imply $(v)$.  Setting  $y=\th^{-r}$ and $c=-\th^{2(j-i+n)},$
the equation $(ii')$ becomes
\begin{eqnarray}\label{Eq4} \th b^2+cy^{10}=-1.  \end{eqnarray}

Completely similar to Case 1, we get $d(\a^{t^{n}S}, \a^{t^{j}\ell S})\ge 2.$

\vskip 3mm
{\it Claim 3: $d(\a^{t^{n}\ell S}, \a^{t^{j}\ell S})\ge 2.$}
\vskip 3mm
The neighborhood of $\a^{t^{n}\ell}$ is
$$\begin{array}{lcl}
Y(i)_1(\a^{t^{n}\ell})&=&\{\a^{t^i\ell ht^{n}\ell}\di h\in H\}=\{\a^{t^{i-n}\ell h_1\ell}\di h_1\in H\}\\
&=&\{\a^g\di g=\overline{\left[\begin{array}{cc} \th^{i-n}&0\\x&\th^{-(i-n)}\end{array}\right]}, x\in
\FF_k\}.\\
\end{array}$$
Then $d(\a^{t^{n}\ell S}, \a^{t^{j}\ell S})=|Y(i)_1(\a^{t^{n}\ell})\cap \a^{t^{j}\ell S}|$, which is
the number of solutions $s(a,b)$  of the equation
\begin{eqnarray}\label{Eq5}
H\overline{\left[\begin{array}{cc} \th^{i-n}&0\\x&\th^{-(i-n)}\end{array}\right]}=Ht^j\ell s(a,b),
\end{eqnarray}
that is
$$\overline{\left[\begin{array}{cc} \th^{5r} &x_1\\ 0 &\th^{-5r} \end{array}\right]}
  \overline{\left[\begin{array}{cc} \th^{i-n}&0\\x&\th^{-(i-n)}\end{array}\right]}=t^j\ell s(a,b),$$
for some $r$ and $x_1$,  that is
$$\left[\begin{array}{cc}
\th^{5r+i-n}+xx_1&x_1\th^{-(i-n)}\\x\th^{-5r}&\th^{-(5r+i-n)}\end{array}\right]
=\pm \left[ \begin{array}{cc} -b\th^{j+1}&-a\th^{j}\\a\th^{-j}&b\th^{-j}\end{array}\right].$$
So Eq(\ref{Eq5}) holds if and only if
$$\left\{\begin{array}{lcll}
\th^{5r+i-n}+xx_1 &=&-b\th^{j+1}, \, \hskip 2cm &(i)\\
x_1\th^{-(i-n)}   &=&-a\th^{j},\, &(ii) \\
x\th^{-5r}          &=&a\th^{-j},\, &(iii) \\
\th^{-(5r+i-n)}   &=&b\th^{-j},\, &(iv)\\
 a^2-b^2\th&=&1,  &(v)\end{array}\right.  $$
that is
$$\left\{\begin{array}{lcll}
\th^{5r+i-n}+xx_1 &=&-b\th^{j+1}, \, \hskip 2cm &(i)\\
x_1   &=&-a\th^{j+i-n},\, &(ii') \\
x      &=&a\th^{5r-j},\, &(iii') \\
b   &=&\th^{j-(5r+i-n)},\, &(iv')\\
 a^2-b^2\th&=&1.  &(v)\end{array}\right.  $$
Inserting $(ii')$, $(iii')$ and $(iv')$ to $(i)$, we obtain
$$a^2\th^{5r+i-n}-\th^{5r+i-n}=\th^{-5r+2j-i+n+1}. \hskip 1cm (i')$$
Note that $(i')$ and $(iv')$ imply $(v)$. Setting  $y=\th^{-r}$ and $c=-\th^{2(j-i+n)-1},$ the
equation $(i')$ becomes $a^2+cy^{10}=1,$ that is Eq(\ref{Eq2}).

Completely similar to Case 1, we get $d(\a^{t^{n}\ell S}, \a^{t^{j}\ell S})\ge 2.$

\vskip 3mm
In summary, Eq(\ref{Eq0}) is  true, which  implies that  the basic orbital graph $Y(i)$ contains a Hamilton cycle.
\qqed

\vskip 3mm
Now we are ready to prove our main theorem.
\vskip 3mm
\f {\bf Proof of Theorem~\ref{the:main}:}
Let $X$ be a graph whose automorphism group contains  a vertex-transitive subgroup  $G:=\PSL(2,s^m)$
where $s$ is a prime,  having a point stabilizer $H:=\ZZ_s^m\rtimes\ZZ_{\frac{s^m-1}{10}}$,
where $s^m+1=2p$ for a  prime $p$.
Then $X$ is an orbital graph $X(G, \OO)$ of $G$ relative to $H$.
Since $G$ has five single point suborbits and  five suborbits of length $k$ which are self-paired, our
$\OO$ contains at least one suborbit $\a^{t^i\ell H}$, otherwise, $X$ is disconnected.
Therefore, $X$ contains a subgraph $Y(i)$, meaning $E(X)$ contains the subset $E(Y(i))$ but $V(X)=V(Y(i))$.
By Lemma~\ref{0}, $Y(i)$ contains a Hamilton cycle, implying that $X$ contains a Hamilton cycle.
\qqed

\end{document}